\documentclass[preprint,12pt]{elsarticle}

\usepackage{amssymb}
\usepackage{amsmath}

\journal{arXiv.org}

\newtheorem{thm}{Theorem}
\newtheorem{lem}[thm]{Lemma}
\newdefinition{rmk}{Remark}
\newproof{pf}{Proof}
\newproof{pot}{Proof of Theorem \ref{thm2}}

\begin{document}

\begin{frontmatter}



\title{Multipoint secant and interpolation methods \\with nonmonotone line search for solving \\systems of nonlinear equations}


\author[label1]{Oleg Burdakov\corref{cor1}}
\ead{oleg.burdakov@liu.se}

\author[label2]{Ahmad Kamandi}
\ead{ahmadkamandi@mazust.ac.ir}

\cortext[cor1]{Corresponding author}

\address[label1]{Department of Mathematics, Link\"oping University, SE-58183 Link\"oping, Sweden}
\address[label2]{Department of Mathematics, University of Science and Technology of Mazandaran, Behshar, Iran}

\begin{abstract}
Multipoint secant and interpolation methods are effective tools for solving systems of nonlinear equations. They use quasi-Newton updates for approximating the Jacobian matrix. Owing to their ability to more completely utilize the information about the Jacobian matrix gathered at the previous iterations, these methods are especially efficient in the case of expensive functions. They are known to be local and superlinearly convergent. We combine these methods with the nonmonotone line search proposed by Li and Fukushima (2000), and study global and superlinear convergence of this combination. Results of numerical experiments are presented. They indicate that the multipoint secant and interpolation methods tend to be more robust and efficient than Broyden's method globalized in the same way.
\end{abstract}


\end{frontmatter}


\section{Introduction}
\label{sec:intr}
Consider the problem of solving a system of simultaneous nonlinear equations
\begin{equation}\label{system} 
F(x)=0,
\end{equation}
where the mapping $F:\mathbb{R}^{n}\rightarrow \mathbb{R}^{n}$ is assumed to be continuously differentiable. Numerical methods aimed at iteratively solving this problem are discussed in \cite{or-70,ds-96,nw-06}. We focus here on those which generate iterates by the formula
\begin{equation}\label{iter} 
x_{k+1} = x_k + \lambda_k p_k, \quad k=0,1, \ldots \ ,
\end{equation}
where the vector $p_k\in \mathbb{R}^{n}$ is a search direction, and the scalar $\lambda_k$ is a step length. 
Denote $F_k = F(x_k)$ and $F'_k = F'(x_k)$. 
In the Newton-type methods, the search direction has the form
$$
p_k = -B_k^{-1}F_k.
$$
Here the matrix $B_k \in \mathbb{R}^{n \times n}$ is either the Jacobian $F'_k$ (Newton's method) or some approximation to it (quasi-Newton methods). For quasi-Newton methods, we consider Broyden's method \cite{b-65}, multipoint secant methods \cite{gs-78,b-83,b-86} and interpolation methods \cite{bf-94,b-97}.

Newton's method is known to attain a local quadratic rate of convergence, when $\lambda_k = 1$ for all $k$.
The quasi-Newton methods do not require computation of any derivatives, and their local rate of convergence is superlinear.

The Newton search direction $p_k^N = - (F'_k)^{-1}F_k$ is a descent direction for $\|F(x)\|$ in any norm. Moreover, as it was shown in \cite{b-80,b-95}, there exists a directional derivative of $\|F(x)\|$ calculated by the formula:
$$
\|F(x_k + \lambda p_k^N)\|^{\prime}_{\lambda = +0} = -\|F_k\|,
$$
which is valid for any norm, even if $\|F(x)\|$ is not differentiable in $x_k$. This property of the Newton search direction provides the basis for constructing various backtracking line search strategies \cite{ds-96,nw-06,b-80} aimed at making Newton's method globally convergent. An important feature of such strategies is that $\lambda_k = 1$ is accepted for all sufficiently large $k$, which allows them to retain the high local convergence rate of the Newton method.

In contrast to Newton's method, the search directions generated by the quasi-Newton methods are not guaranteed to be descent directions for $\|F(x)\|$. This complicates the globalization of the latter methods. 

The earliest line search strategy designed for globalizing Broyden's method is due to Griewank \cite{g-86}. Its drawback, as indicated in \cite{lf-00}, is related to the case when $p_k$ is orthogonal, or close to orthogonal, to the $\nabla \|F(x_k)\|^2$. Here and later, $\|\cdot \|$ stands for the Euclidean vector norm and the induced matrix norm. The Frobenius matrix norm will be denoted by $\|\cdot \|_F$.

Li and Fukushima \cite{lf-00} developed a new backtracking line search for Broden's method and proved its global superlinear convergence. In this line search, the function $\|F_k\|$ may not monotonically decrease with $k$. Its important feature is that it is free of the aforementioned drawback of the line search proposed in \cite{g-86}.

The purpose of this paper is to extend the Li-Fukushima line search to the case of the multipoint secant and interpolation methods, theoretically study their global convergence and also explore their practical behavior in numerical experiments. We are also aimed at demonstrating a higher efficiency of these methods as compared with Broyden's method in the case of expensive function evaluations.

The paper is organized as follows. In the next section, we describe the multipoint secant and interpolation methods and discuss their properties. A combination of these methods with the Li-Fukushima line search is presented in Section~\ref{sec:LineSearch}. In Section~\ref{sec:Convergence}, we show a global and superlinear convergence of this combination.    
Results of numerical experiments are reported and discussed in Section~\ref{sec:NumExp}.
Finally, some conclusions are included in the last section of the paper.

\section{Quasi-Newton updates}
\label{sec:QN}
The class of quasi-Newton updates that we consider here has the form
\begin{equation}\label{QNupdates} 
B_{k+1}=B_k+\frac{(y_k-B_ks_k)c_k^T}{s_k^Tc_k},
\end{equation}
where $s_k=x_{k+1}-x_k$, $y_k=F_{k+1}-F_k$, and $c_k \in \mathbb{R}^{n}$ is a parameter.

One of the most popular quasi-Newton method of solving \eqref{system} is due to Broyden \cite{b-65}. It corresponds to the choice $c_k=s_k$ and satisfies the, so-called, \emph{secant equation}:
\begin{equation}\label{sec_eq}  
B_{k+1}s_k=y_k.
\end{equation}
It indicates that $B_{k+1}$ provides an approximation of the Jacobian matrix along the direction $s_k$. Though such an approximation is provided by $B_k$ along $s_{k-1}$, it is not guaranteed that $B_{k+1}$ retains this property because, in general, $B_{k+1}s_{k-1} \neq y_{k-1}$.

Gay and Schnabel \cite{gs-78} proposed a quasi-Newton updating formula of the form \eqref{QNupdates} with the aim to preserve the secant equations satisfied at some previous iterations. The resulting Jacobian approximation satisfies the following \emph{multipoint secant equations}:
\begin{equation}\label{multi-secant} 
B_{k+1}s_i=y_i, \quad  \forall i\in T_{k+1},
\end{equation}
where $T_{k+1}=\{i: m_k\leq i  \leq k\}$ and $ 0\leq m_k  \leq k$. To guarantee this, the parameter in \eqref{QNupdates} is calculated by the formula
\begin{equation}\label{c_k}
	c_k = s_k - P_ks_k,
\end{equation}
where $P_k \in \mathbb{R}^{n \times n}$ is an orthogonal projector on the subspace generated by the vectors 
$s_{m_k}, s_{m_k+1},\ldots ,s_{k-1}$, and $P_k$ vanishes when $m_k = k$. To ensure a local superlinear convergence and stable approximation of the Jacobian, it is required in \cite{gs-78} that there exists $\bar{\sigma} \in (0,1)$ such that 
\begin{equation}\label{sin}
	\|c_k\| \ge \bar{\sigma} \|s_k\|, \quad \forall k\ge 0.
\end{equation}
To meet this requirement, $m_k$ is chosen as follows. If the trial choice of $m_k=m_{k-1}$ fails to satisfy \eqref{sin}, the vectors $s_{m_{k-1}}, \ldots ,s_{k}$ are considered as close to linear dependent, and then a \emph{restart} is performed by setting $m_k=k$, or equivalently, $T_{k+1} = \{k\}$. Otherwise, the trial choice is accepted, in which case the set $T_{k+1}$ is obtained by adding $\{k\}$ to $T_k$.

In what follows, we say, for a given $\sigma \in (0,1)$, that non-zero vectors $v_i\in \mathbb{R}^{n}$, $i=1,2, \ldots ,m$, are \emph{$\sigma$-safely linearly independent} if the inequality
\begin{equation}\label{lin_ind}
	\det \left(
	\left[\frac{v_1}{\|v_1\|}, \ldots \frac{v_m}{\|v_m\|}\right]^T
	\left[\frac{v_1}{\|v_1\|}, \ldots \frac{v_m}{\|v_m\|}\right]
	\right) \ge \sigma^2
\end{equation}
holds. Here the ordering of the vectors is not essential. Note that, for each $k$ in the Gay-Schnabel method, the vectors $\{s_i\}_{i\in T_k}$ are $\sigma$-safely linearly independent, where $\sigma$ depends only on $\bar{\sigma}$ and $n$.

It should be mentioned that, in the case of restart, the multipoint secant equations \eqref{multi-secant} are reduced to the single secant equation \eqref{sec_eq}, which means that the collected information about the Jacobian is partially lost. The quasi-Newton methods proposed in \cite{b-83,b-86} are aimed at avoiding restarts.  In these methods, the vectors $\{s_i\}_{i\in T_k}$ are also $\sigma$-safely linearly independent. Instead of setting $T_{k+1} = \{k\}$, when the vectors $\{s_i\}_{i\in T_k \cup \{k\}}$ do not meet this requirement, the set $T_{k+1}$ is composed of those indices in $T_k$ which, along with the index $k$, ensure that $\{s_i\}_{i\in T_{k+1}}$ are $\sigma$-safely linearly independent. Since the way of doing this is not unique, a preference may be given, for instance, to the most recent iterations in $T_k$ because they carry the most fresh information about the Jacobian. The Jacobian approximation is updated by formula \eqref{QNupdates} with $c_k$ computed in accordance with \eqref{c_k}, where $P_k$ is the orthogonal projector on the subspace generated by the vectors $\{s_i\}_{i\in T_{k+1} \setminus \{k\}}$. The methods in \cite{b-83,b-86} are superlinearly convergent.

For describing the interpolation methods, we need the following definition. For a given $\sigma \in (0,1)$, we say that points $x_i\in \mathbb{R}^{n}$, $i=0,1, \ldots ,m$, are in \emph{$\sigma$-stable general position} if there exist vectors $\{\Delta x_j\}_{j=1}^{m}$ of the form $x_{\mu_j}-x_{\nu_j}$, $0 \leq \mu_j, \nu_j \leq m$ such that they are $\sigma$-safely linearly independent, which means that the inequality
\begin{equation}\label{gen_pos_1} 
	\det (\Delta X^T \Delta X)
	\ge \sigma^2
\end{equation}
holds for the matrix
$$
\Delta X = \left[\frac{\Delta x_1}{\|\Delta x_1\|}, \ldots , \frac{\Delta x_m}{\|\Delta x_m\|}\right].
$$
Here the ordering of the vectors is not essential, whereas a proper choice of such vectors does. The latter is equivalent to choosing a most linearly independent set of $m$ vectors of the form $x_{p_j}-x_{q_j}$ which constitute a basis for the linear manifold generated by the points $\{x_i\}_{i=0}^m$. In \cite{b-97}, an effective algorithm for finding vectors, which minimizes the value of the left-hand side in \eqref{gen_pos_1}, was introduced. It is based on a reduction of this minimization problem to a minimum spanning tree problem formulated for a graph whose nodes and edges correspond, respectively, to the points and all the vectors connecting the points. Each edge cost is equal to the length of the respective vector. It is also shown in \cite{b-97} how to effectively update the minimal value of the determinant when one point is removed from or added to the set.

As it was pointed out in \cite{bf-94,b-97}, when search directions are close to be linearly dependent, the corresponding iterates still may be in a stable general position, which provides a stable Jacobian approximation. In such cases, instead of discarding some information about the Jacobian provided by the pairs $(s_i, y_i)$, the quasi-Newton methods introduced in \cite{bf-94,b-97} make use of this kind of information provided by the pairs $(x_i, F_i)$. At iteration $k$, they construct an interpolating linear model $L_k(x)=F_k+B_k(x-x_k)$ such that
\begin{equation}\label{interpolation} 
L_k(x_i)=F_i, \quad \forall i\in I_k,
\end{equation}
where $I_k$ is a set of indices with the property that $\{k,k-1\}\subseteq I_k \subseteq \{k,k-1,...,k-n\}$. Then the solution to the system of linear equations $L_k(x)=0$ yields the new iterate $x_{k+1}$. The Jacobian approximation is updated by formula \eqref{QNupdates}, in which
$$
c_k = x_{k+1} - P_kx_{k+1},
$$
where $P_k$ is the orthogonal projector on the linear manifold generated by the points $\{x_i\}_{i \in I_{k+1} \setminus \{k+1\}}$. The interpolation property is maintained by virtue of including in $I_{k+1}$ elements $\{k,k+1\}$ and some elements of the set $I_{k}$. The main requirement, which ensures a stable Jacobian approximation and superlinear convergence, is that the iterates $\{x_i\}_{i \in I_{k+1}}$ are in the $\sigma$-stable general position. Since the way of choosing indices of $I_{k}$ for including in $I_{k+1}$ is not unique, it is desirable to make a priority for the most recent iterates.

The only difference between the quasi-Newton methods considered here is in their way of computing the vector $c_k$. For Broyden's method, it is the least expensive, whereas the multipoint secant and interpolation methods require, as one can see in Section~\ref{sec:NumExp}, far less number of function evaluations. Therefore, the latter quasi-Newton methods are more suitable for solving problems, in which one function evaluation is more expensive than the computation of $c_k$.

The computational cost of each iteration in the considered quasi-Newton methods depends on the number of couples $(s_i,y_i)$ or $(x_i,F_i)$ that are involved in calculating $c_k$. Therefore, in some problems, especially those of large scale, it is reasonable to limit the number of stored couples by limiting the depth of memory. This can be done by introducing a parameter $m\le n$ which prevents from using the couples with $i<k-m$. Note that, Broyden's method is a special case of the multipoint secant and interpolation methods for $m=0$ and $m=1$, respectively.

In the considered quasi-Newton methods, the matrix $B_{k+1}$, like in Broyden's method, results from a least-change correction to $B_{k}$ in the Frobenius norm over all matrices that satisfy the corresponding secant or interpolation conditions. A related property, which is common to these methods, is that the vector $c_k$ in \eqref{QNupdates} is such that $c_k^Ts_k = \|c_k\|^2$. 

\section{Quasi-Newton algorithms with Li-Fukushima line search}
\label{sec:LineSearch}
In this section, we present the Li-Fukushima line search \cite{lf-00} adapted to the class of the quasi-Newton updates considered above. The matrix $B_{k+1}$ is normally nonsingular. If not, it is computed by the modified updating formula
\begin{equation}\label{nonsingularB}
	B_{k+1}=B_k+\theta_k \frac{(y_k-B_ks_k)c_k^T}{\|c_k\|^2}.
\end{equation}
Here $\theta_k \in [1 - \bar{\theta}, 1 + \bar{\theta}]$ is chosen so that $B_{k+1}$ is nonsingular, where the parameter $\bar{\theta} \in (0,1)$. 

It should be noted that the theoretical analysis of formula \eqref{nonsingularB} conducted in \cite{lf-00} for $c_k=s_k$ points to the interesting fact that Broyden's methods retains its superlinear convergence, even if to fix $\theta_k \in [1 - \bar{\theta}, 1 + \bar{\theta}]$ for all $k$, provided that the resulting $B_{k+1}$ is nonsingular at all iterations. In this case, the secant equation \eqref{sec_eq} is not necessarily satisfied.

The Li-Fukushima line search involves a positive sequence $\{\eta_k\}$ such that
\begin{equation}\label{sum_eta} 
\sum_{k=0}^{\infty}\eta_k < \infty.
\end{equation}
The line search consists in finding a step length $\lambda$ which satisfies the inequality
\begin{equation}\label{nonmonot_ls}  
\|F(x_k+\lambda p_k)\|\leq \|F_k\|-\sigma_1\|\lambda p_k\|^2+\eta_k\|F_k\|,
\end{equation}
where $\sigma_1>0$ is a given parameter. This inequality is obviously satisfied for all sufficiently small values of $\lambda >0$ because, as $\lambda$ goes to zero, the left-hand and right-hand sides of \eqref{nonmonot_ls} tend to $\|F_k\|$ and $(1+\eta_k)\|F_k\|$, respectively.

A step length which satisfies \eqref{nonmonot_ls} can be produced by the following backtracking procedure.
\begin{tabbing}
\rule{\textwidth}{1px}\\
{\bf Algorithm~1} \ Backtracking procedure.\\
\rule{\textwidth}{0.5px}\\
\hspace{0.5cm}\=\hspace{0.5cm}\=\hspace{0.5cm}\=\hspace{0.5cm}\=\kill
\textbf{Given:} $\sigma_1>0$, $\eta_k>0$ and $\beta \in (0,1)$\\
Set $\lambda \leftarrow 1$\\
\textbf{repeat} until \eqref{nonmonot_ls} is satisfied\\
\> $\lambda \leftarrow \beta \lambda$\\
\textbf{end (repeat)}\\
\textbf{return} $\lambda_k = \lambda$\\
\rule{\textwidth}{0.5px}
\end{tabbing}

Note that the Li-Fukushima line search is nonmonotone because the monotonic decrease $\|F_{k+1}\| < \|F_k\|$ may be violated at some iterations. Since $\eta_k \rightarrow 0$, the size of possible violation of monotonicity vanishes. This line search is extended below to the case of the quasi-Newton methods considered in the present paper.
\begin{tabbing}
\rule{\textwidth}{1px}\\
{\bf Algorithm~2} \ Quasi-Newton methods with Li-Fukushima line search.\\
\rule{\textwidth}{0.5px}\\
\hspace{0.5cm}\=\hspace{0.5cm}\=\hspace{0.5cm}\=\hspace{0.5cm}\=\kill
\textbf{Given:} initial point $x_0 \in \mathbb{R}^{n}$, nonsingular matrix $B_0 \in \mathbb{R}^{n \times n}$, positive sca-\\ lars $\sigma_1, \sigma_2 >0$, \ $\beta, \sigma, \rho, \bar{\theta} \in (0,1)$, and positive sequence $\{\eta_k\}$ satisfying \eqref{sum_eta}.\\
\textbf{for} $k=0,1,2, \ldots$ \textbf{do}\\
\> \textbf{if} $F_k=0$ \textbf{then} stop. \\
\> Find $p_k$ that solves $B_kp + F_k =0$.\\
\> \textbf{if} $\|F(x_k+p_k)\|\leq \rho \|F_k\|-\sigma_2\|p_k\|^2$ \textbf{then} set $\lambda_k \leftarrow 1$ \textbf{else} \\
\> \> use Algorithm~1 for finding $\lambda_k$.\\
\> \textbf{end (if)}\\
\> Compute $c_k \in \mathbb{R}^{n}$ in accordance with the chosen quasi-Newton method.\\
\> Compute nonsingular $B_{k+1}$ by properly choosing $\theta_k \in [1 - \bar{\theta}, 1 + \bar{\theta}]$ in \eqref{nonsingularB}. 
\\
\textbf{end (for)}
\\
\rule{\textwidth}{0.5px}
\end{tabbing}

As it was mentioned above, in Broyden's method, $c_k=s_k$. We present now generic algorithms of computing $c_k$
for the multipoint secant and interpolation methods separately. 

The multipoint secant methods \cite{gs-78,b-83,b-86} start with the set $T_0 = \emptyset$, and then they proceed in accordance with the following algorithm.
\begin{tabbing}
\rule{\textwidth}{1px}\\
{\bf Algorithm~3} \ Computing $c_k$ for the multipoint secant methods.\\
\rule{\textwidth}{0.5px}\\
\hspace{0.5cm}\=\hspace{0.5cm}\=\hspace{0.5cm}\=\hspace{0.5cm}\=\kill
\textbf{Given:} $\sigma \in (0,1)$, $s_k$ and $\{s_i\}_{i\in T_k}$.\\
Set $T_k \leftarrow T_k \setminus \{k-n\}$.\\
Find $T_{k+1} \subseteq T_k \cup \{k\}$ such that $\{k\} \subseteq T_{k+1}$, and $\{s_i\}_{i\in T_{k+1}}$ are $\sigma$-safely \\ linearly independent.\\
Set $c_k \leftarrow s_k - P_ks_k$, where $P_k$ is the orthogonal projector onto the subspace \\generated by $\{s_i\}_{i\in T_{k+1}\setminus\{k\}}$.
\\
\textbf{return} $c_k$ and $\{s_i\}_{i\in T_{k+1}}$.\\
\rule{\textwidth}{0.5px}
\end{tabbing}

In the interpolation methods \cite{bf-94,b-97}, the initial set $I_0 = \{0\}$. They are based on the following algorithm.
\begin{tabbing}
\rule{\textwidth}{1px}\\
{\bf Algorithm~4} \ Computing $c_k$ for the interpolation methods.\\
\rule{\textwidth}{0.5px}\\
\hspace{0.5cm}\=\hspace{0.5cm}\=\hspace{0.5cm}\=\hspace{0.5cm}\=\kill
\textbf{Given:} $\sigma \in (0,1)$, $x_{k+1}$ and $\{x_i\}_{i\in I_k}$.\\
Set $I_k \leftarrow I_k \setminus \{k-n\}$.\\
Find $I_{k+1} \subseteq I_k \cup \{k+1\}$ such that $\{k, k+1\} \subseteq I_{k+1}$, and $\{x_i\}_{i\in I_{k+1}}$ are in \\$\sigma$-stable general position.\\
Set $c_k \leftarrow x_{k+1} - x_{k+1}^{\perp}$, where $x_{k+1}^{\perp}$ is the orthogonal projection of the point \\$x_{k+1}$ onto 
the linear manifold generated by $\{x_i\}_{i\in I_{k+1}\setminus\{k+1\}}$.
\\
\textbf{return} $c_k$ and $\{x_i\}_{i\in I_{k+1}}$.\\
\rule{\textwidth}{0.5px}
\end{tabbing}

Algorithms~3 and 4 pose certain restrictions on choosing the sets $T_{k+1}$ and $I_{k+1}$, respectively. However, they also admit some freedom in choosing the sets. In this sense, each of these algorithms represents a class of methods. Specific choices of the sets and implementation issues are discussed in Section~\ref{sec:NumExp}. Note that $T_{k+1} = \{k\}$ and $I_{k+1} = \{k, k+1\}$ are valid choices which result in $c_k = s_k$. This means that Broyden's method is a special case of the two classes. Therefore, the convergence analysis presented in the next section can be viewed as an extension of the results in \cite{lf-00}. It should be emphasized that the extension is not straightforward, because it requires establishing some nontrivial features of the multipoint secant and interpolation methods.

\section{Convergence analysis}
\label{sec:Convergence}
To study the convergence of the quasi-Newton methods with Li-Fukushima line search, we will use the next three lemmas proved in \cite{lf-00}. They do not depend on the way of generating the search directions $p_k$.

\begin{lem}\label{lem1}
The sequence $\{x_k\}$ generated by Algorithm~2 is contained in the
set
\begin{equation}\label{Omega} 
\Omega=\{x\in \mathbb{R}^n : \ \|F(x)\|\leq e^{\eta}\|F_0\|\},
\end{equation}
where 
$$
\eta = \sum_{k=0}^{\infty}\eta_k.
$$
\end{lem}

\begin{lem}\label{lem2}
Let the level set $\Omega$ be bounded and $\{x_k\}$ be generated by Algorithm~2. Then
\begin{equation}\label{e:19} 
\sum_{k=0}^{\infty}\|s_k\|^2 < \infty.
\end{equation}
\end{lem}

\begin{lem}\label{a_b_xi}
Let $\{a_k\}$, $\{b_k\}$ and $\{\xi_k\}$ be positive sequences satisfying
$$ 
a_{k+1}^2 \le (a_k + b_k)^2 - \alpha \xi_k^2, \quad k = 0,  1, \ldots ,
$$ 
where $\alpha$ is a constant. Then 
\begin{equation}\label{sum_b_squared}
	\sum_{k=0}^{\infty} b_k^2 < \infty \quad \Rightarrow \quad 
  \lim_{k \rightarrow \infty} \frac{1}{k} \sum_{i=0}^{k-1} \xi_i^2 = 0,
\end{equation}
and
\begin{equation}\label{sum_b}
	\sum_{k=0}^{\infty} b_k < \infty \quad \Rightarrow \quad 
	\sum_{k=0}^{\infty} \xi_k^2 \le \infty .
\end{equation}
\end{lem}

The further convergence analysis requires the following assumptions.

\begin{description}
	\item[A1.] 
The level set $\Omega$ defined by \eqref{Omega} is bounded.
	\item[A2.] 
The Jacobian $F'(x)$ is Lipschitz continuous on the convex hall of $\Omega$, i.e., there exists a positive constant $L$ such that
$$
\|F'(x)-F'(y)\| \leq L\|x-y\|, \qquad \forall x,y \in  \textrm{Conv}(\Omega).
$$
	\item[A3.] 
$F'(x)$ is nonsingular for every $x \in \Omega$.
\end{description}

The set $\Omega$ in A2 and A3 is not necessarily assumed to be the level set defined by \eqref{Omega}, unless these assumptions are combined with A1 in one and the same assertion.

We begin with establishing global convergence result for the interpolation methods represented by Algorithm~2, in which $c_k$ is produced by Algorithm~4. By construction, the interpolation points $\{x_i\}_{i\in I_{k+1}}$ are in $\sigma$-stable general position. This means that there exist $\ell_{k+1} = |I_{k+1}|-1$ vectors, $\{\Delta x_j\}_{j=1}^{\ell_{k+1}}$, such that, first, they are of the form $\Delta x_j = x_{\mu_j} - x_{\nu_j}$, where $\mu_j, \nu_j \in I_{k+1}$, and second, inequality \eqref{gen_pos_1} holds for the corresponding matrix
$$
\Delta X = 
\left[\frac{\Delta x_1}{\|\Delta x_1\|}, \ldots , \frac{\Delta x_{\ell_{k+1}}}{\|\Delta x_{\ell_{k+1}}\|}\right].
$$
Let $\Delta X_{\perp} \in \mathbb{R}^{n \times (n-\ell_{k+1})}$ be an orthonormal matrix such that $\Delta X_{\perp}^T \Delta X = 0$. Denote
\begin{equation}\label{e:21} 
A_{k+1} = \sum_{j=1}^{\ell_{k+1}} \frac{\Delta F_j u_j^T}{\Delta x_j^Tu_j} 
+ F_{k+1}^{\prime} \Delta X_{\perp}\Delta X_{\perp}^T,
\end{equation}
where $F_j = F_{\mu_j} - F_{\nu_j}$, $u_j=\Delta x_j-P_j\Delta x_j$, and $P_j$ is the orthogonal projector onto the subspace generated by all the vectors $\Delta x_1, \ldots , \Delta x_{\ell_{k+1}}$, except the vector $\Delta x_j$. It follows from \eqref{e:21} that
\begin{equation}\label{e:22} 
A_{k+1} \Delta x_j = \Delta F_j, \qquad j=1, \ldots , \ell_{k+1}.
\end{equation}
Consequently, 
\begin{equation}\label{e:23}
A_{k+1}(x_i - x_j) = F_i - F_j, \qquad \forall i,j\in I_{k+1}.
\end{equation}

The next result establishes a key property of the matrix $A_{k+1}$. In its formulation, we disregard the way in which the iterates are generated.
The property of $A_{k+1}$ will be used for showing global convergence of the interpolation methods.

\begin{lem}\label{lem3} 
Let points $\{x_i\}_{i\in I_{k+1}} \subseteq \{x_i\}_{i=k-n+1}^{k+1}$ be in $\sigma$-stable general position. Suppose that assumption A2 holds for the set 
$$
\Omega = \{x_i\}_{i=k-n}^{k+1}.
$$
Then 
\begin{equation}\label{e:24} 
\|A_{k+1} - F_{k+1}^{\prime}\| \leq \frac{L \sqrt{n}}{\sigma} \sum_{i=k-n+1}^{k}\|s_i\|.
\end{equation}
If, in addition, points $\{x_i\}_{i\in I_{k}} \subseteq \{x_i\}_{i=k-n}^{k}$ are also in $\sigma$-stable general position and belong to $\Omega$, then
\begin{equation}\label{deltaA}
	\|A_{k+1} - A_k\| \le \frac{3L \sqrt{n}}{\sigma} \sum_{i=k-n}^{k}\|s_i\|.
\end{equation}
\end{lem}
\begin{pf}
Consider the matrix $Q = [\Delta X \ \ \Delta X_{\perp}]$. It can be easily shown that 
\begin{equation}\label{invQ}
	\|Q^{-1}\| \le 1/ \sigma.
\end{equation}
Indeed, the upper bound in \eqref{invQ} is related to the smallest eigenvalue of the matrix $Q^TQ$, which is a block-diagonal matrix, whose two blocks are $\Delta X^T \Delta X$ and the identity matrix of the proper size. From the fact that the smallest eigenvalue of the first block is bounded below by $\sigma^2$, we get \eqref{invQ}.

Note that $(A_{k+1} - F_{k+1}^{\prime}) \Delta X_{\perp} = 0$, and
$$
(A_{k+1} - F_{k+1}^{\prime}) \Delta X = \left[
\frac{\Delta F_1 - F_{k+1}^{\prime}\Delta x_1}{\|\Delta x_1\|}, \ldots , 
\frac{\Delta F_{\ell_{k+1}} - F_{k+1}^{\prime}\Delta x_{\ell_{k+1}}}{\|\Delta x_{\ell_{k+1}}\|}
\right] ,
$$
For the columns of this matrix, \cite[Theorem~3.2.5]{or-70} and assumption A2 give 
\begin{equation}\label{Delta_F-F'Delta_x}
\frac{\|\Delta F_j - F_{k+1}^{\prime}\Delta x_j\|}{\|\Delta x_j\|}
\le L \max \{\|x_{\mu_j} - x_{k+1}\|, \|x_{\nu_j} - x_{k+1}\|\}, \quad j = 1, \ldots , \ell_{k+1}.
\end{equation}

Then, using a matrix norm equivalence \cite[Theorem 3.3]{kl-95}, \eqref{invQ} and \eqref{Delta_F-F'Delta_x}, we get
\begin{align}
\|A_{k+1} - F_{k+1}^{\prime}\| &= \|(A_{k+1} - F_{k+1}^{\prime})Q Q^{-1}\| 
\le \|(A_{k+1} - F_{k+1}^{\prime})\Delta X\| \|Q^{-1}\|
\nonumber \\
&\le \frac{L \sqrt{n}}{\sigma} \max_{1 \le j \le \ell_{k+1}} \{\|x_{\mu_j} - x_{k+1}\|, \|x_{\nu_j} - x_{k+1}\|\}.
\nonumber
\end{align}
From this inequality one can easily conclude that \eqref{e:24} holds.

Observe that 
$$
\|A_{k+1} - A_k\| \le \|A_{k+1} - F_{k+1}^{\prime}\| + \|F_{k+1}^{\prime} - F_{k}^{\prime}\|
+ \|A_{k} - F_{k}^{\prime}\|.
$$
This inequality along with assumption A2 and inequality \eqref{e:24} show that \eqref{deltaA} holds, so our
proof is complete.
\qed \end{pf}

Consider the interpolation property \eqref{interpolation}. It implies that
\begin{equation}\label{general_secant}
	B_{k}(x_i - x_j) = F_i - F_j, \qquad \forall i,j\in I_{k}.
\end{equation}
Similar relations are established for $A_{k+1}$ in \eqref{e:23}. They hold in particular for all $i,j \in I_{k+1}\setminus \{k+1\}$. By construction, $I_{k+1}\setminus \{k+1\} \subseteq I_k$. Then, combining \eqref{e:23} and \eqref{general_secant}, we get the relation
$$
B_{k}(x_i - x_j) = A_{k+1}(x_i - x_j), \quad \forall i,j \in I_{k+1}\setminus \{k+1\}.
$$
Hence,
$$
B_{k}(x' - x'') = A_{k+1}(x' - x''), \quad \forall x', x'' \in \mathcal{L},
$$
where $\mathcal{L}$ is the linear manifold generated by the points $\{x_i\}_{i \in I_{k+1}\setminus \{k+1\}}$. It is easy to see that this relation yields $B_{k}(c_k - s_k) = A_{k+1}(c_k - s_k)$, or equivalently,
\begin{equation}\label{B-A}
(B_{k} - A_{k+1})c_k = (B_{k} - A_{k+1})s_k.
\end{equation}
Indeed, recall that $s_k = x_{k+1} - x_k$ and $c_k = x_{k+1} - x_{k+1}^{\perp}$, which means that $c_k - s_k = x_k - x_{k+1}^{\perp}$, where $x_k, x_{k+1}^{\perp} \in \mathcal{L}$.

Note that the equation $A_{k+1}s_k=y_k$ is a special case of \eqref{e:23}. Then the updating formula \eqref{nonsingularB} can be written as
\begin{equation}\label{e:28} 
B_{k+1} = B_k + \theta_k \frac{(A_{k+1} - B_k)s_k c_k^T}{\|c_k\|^2}.
\end{equation}
By analogy with \cite{lf-00}, we define
$$
\xi_k = \frac{\|y_k - B_k s_k\|}{\|c_k\|}.
$$
In the next result, which is similar to \cite[Lemma 2.6]{lf-00}, we study the behaviour of this sequence in the case of $B_k$ generated by the interpolation methods.

\begin{lem}\label{lem4}
Let assumptions A1 and A2 hold, and $\{x_k\}$ be generated by Algorithms~2 and 4. If
\begin{equation}\label{e:30} 
\sum_{k=0}^{\infty}\|s_k\|^2<\infty ,
\end{equation}
then
\begin{equation}\label{e:31} 
\lim_{k \rightarrow\infty}\frac{1}{k}\sum_{i=0}^{k-1}\xi_i^2=0 .
\end{equation}
In particular, there exists a subsequence of $\{\xi_k\}$ which converge to zero. If
\begin{equation}\label{e:32} 
\sum_{k=0}^{\infty}\|s_k\|<\infty ,
\end{equation}
then
\begin{equation}\label{e:33} 
\sum_{k=0}^{\infty}\xi_k^2<\infty .
\end{equation}
In particular, the whole sequence $\{\xi_k\}$ converge to zero.
\end{lem}
\begin{pf}
Denote
$$
a_k=\|B_k-A_k\|_F \quad \textrm{and} \quad b_k=\|A_{k+1}-A_k\|_F.
$$
From \eqref{e:28}, we have
\begin{equation*}\label{e:34}
\begin{split}   
a_{k+1}^2 
&= \left| \left|(B_k - A_{k+1})\left(I - \theta_k \frac{s_kc_k^T}{\|c_k\|^2}\right) \right| \right|_F^2\\
&= \|B_k-A_{k+1}\|_F^2 
- 2 \theta_k \textrm{trace}\left(\frac{(B_k - A_{k+1})s_kc_k^T(B_k - A_{k+1})^T}{\|c_k\|^2}\right)\\
& \hspace{3.12cm} + \theta_k^2 \textrm{trace}\left(\frac{(B_k-A_{k+1})s_ks_k^T(B_k-A_{k+1})^T}{\|c_k\|^2}\right).
\end{split} 
\end{equation*}
Using here \eqref{B-A}, we get 
\begin{equation*}\label{e:35}
a_{k+1}^2 = \|B_k - A_{k+1}\|_F^2 - \theta_k(2 - \theta_k) \frac{\|(B_k - A_{k+1})s_k\|^2}{\|c_k\|^2}.
\end{equation*}
The triangular inequality yields $\|B_k - A_{k+1}\|_F^2 \le (a_k + b_k)^2$. Furthermore, 
$\theta_k(2 - \theta_k) \ge (1 - \bar{\theta}^2) > 0$, because $|\theta_k - 1| \le \bar{\theta}$. Then
\begin{equation*}\label{e:36}
a_{k+1}^2 \le (a_k+b_k)^2 - (1 - \bar{\theta}^2)\xi_k^2.
\end{equation*}
This inequality ensures that the main assumption of Lemma~\ref{a_b_xi} holds. Let condition \eqref{e:30} be satisfied. Then, by Lemma~\ref{lem3} and norm equivalence, the implication \eqref{sum_b_squared} is applicable, which proves \eqref{e:31}. Supposing now that condition \eqref{e:32} is satisfied, we can similarly show that the implication \eqref{sum_b_squared} is applicable, and it yields \eqref{e:33}. This completes the proof.
\qed \end{pf}

It can be easily seen that the results obtained so far for the interpolation methods are also valid in the case of the multipoint secant methods. This can be verified by substituting $s_j = x_{j+1} - x_j$ for $\Delta x_j$ in \eqref{e:21} and also in the subsequent manipulations with $\Delta x_j$. 

We are now in a position to derive convergence results for the multipoint secant and interpolation methods globalized by means of Algorithm~2.

\begin{thm}\label{thm1} Let assumptions A1, A2 and A3 hold. Suppose that the sequence $\{x_k\}$ is generated by Algorithm~2, where the vector $c_k$ is produced by either of Algorithms 3 or 4. Then $\{x_k\}$ converges to the unique solution of \eqref{system}. Moreover, the rate of convergence is superlinear.
\end{thm}
\begin{pf}
We skip the proof of convergence to the unique solution of \eqref{system} because it is entirely similar to that of \cite[Theorem 2.1]{lf-00}. One major difference is that the quantity 
$$
\zeta_k = \frac{\|y_k - B_k s_k\|}{\|s_k\|}.
$$
is used in \cite{lf-00} instead of the $\xi_k$ that is used in the present paper. The relation between the two quantities is the following. The vector $c_k$ generated by Algorithms 3 and 4 is such that $c_k^T s_k = \|c_k\|^2$, that is $\|c_k\| \le \|s_k\|$. Thus, $\zeta_k \le \xi_k$, and therefore, the statements of Lemma~\ref{lem4} refer also to the sequence $\{\zeta_k\}$. This allows us to invoke here \cite[Theorem 2.1]{lf-00}.

We skip the proof of superlinear convergence because it follows the same steps as in \cite[Theorem 2.2]{lf-00}.
\qed \end{pf}

This result shows that the globalized multipoint secant and interpolation methods have the same theoretical properties as Broyden's method. However, as one can see in the next section, the former methods have some practical advantages.

\section{Numerical experiments}
\label{sec:NumExp}
The developed here global convergent quasi-Newton algorithms were implemented in MATLAB. We shall refer to them as
\begin{description}
	\item[QN1:] Broyden's method \cite{b-65},
	\item[QN2:] Gay-Schnabel's multipoint secant method \cite{gs-78},
	\item[QN3:] multipoint secant method \cite{b-83,b-86},
	\item[QN4:] interpolation method \cite{bf-94,b-97}.
\end{description}
Each of them is a special case of Algorithm~2. The difference between them consists in the following specific ways of computing the parameter $c_k$. 
\begin{description}
\item[QN1:] $c_k \leftarrow s_k$. 

\item[QN2:] The parameter $c_k$ is computed by Algorithm~3 as follows.\\
Set $T_{k+1} \leftarrow T_k \cup \{k\}$ and $c_k \leftarrow s_k - P_ks_k$.  \\
\textbf{if} $\|c_k\| \le \sigma \|s_k\|$ \textbf{then} $T_{k+1} \leftarrow \{k\}$ and $c_k \leftarrow s_k$.

\item[QN3:] The parameter $c_k$ is computed by Algorithm~3 as follows.\\
Set $S_k \leftarrow [\ldots,\frac{s_i}{\|s_i\|},\ldots]_{i\in T_{k}\cup \{k\}}$, where the columns are sorted in decreasing order of the indices.\\
Compute $QR$ factorization of $S_k$ so that all diagonal elements of $R$ are non-negative. \\
Compute $d_k = \det (S_k^T S_k) = \prod_{i \in T_k} R_{ii}^2$, where $R_{ii}$ is the diagonal element of $R$ that corresponds to the column $s_i/\|s_i\|$. \\
\textbf{while} $d_k < \sigma^2$ \textbf{do} \\
  \hspace*{4mm} Find $j = \textrm{arg\,min}\{R_{ii}: \ i\in T_k\}$. \\
  \hspace*{4mm} Set $T_{k} \leftarrow T_k \setminus \{j\}$ and
	compute $d_k = \prod_{i \in T_k} R_{ii}^2$ (or, equivalently, set \\
	\hspace*{4mm} $d_k  \leftarrow d_k / R_{jj}^2$ when $R_{jj} \neq 0$).\\
\textbf{end while}\\
Set $T_{k+1} \leftarrow T_k \cup \{k\}$ and $c_k \leftarrow s_k - P_ks_k$.

\item[QN4:] The parameter $c_k$ is computed by Algorithm~4 in which the set $I_{k+1}$ is produced in accordance with \cite[Algorithm~4.1]{b-97}.
\end{description}

Note that in the while-loop of QN3, $QR$ is not computed for any new matrix $S_k$. Since the columns of $S_k$ are of unit length, all diagonal elements of $R$ are such that $R_{ii} \in [0,1]$ with $R_{kk} = 1$. In this connection, it can be easily seen that if to remove any column in $S_k$, then the diagonal elements of the new $R$-factor (if computed) cannot be smaller than the corresponding old ones. Thus, at any step of the while-loop, we have $d_k \le \det (S_k^T S_k)$. Consequently, the vectors $\{s_i\}_{i\in T_{k+1}}$ obtained by QN3 are $\sigma$-safely linearly independent.

In the four algorithms, the stopping criterion was
$$
\|F(x_k)\| \le 10^{-10}\cdot\max\{\|F(x_0)\|,1\}.
$$
The parameters were chosen as $\sigma=0.1$, $\sigma_1=\sigma_2=0.001$, $\rho =0.9$, $\beta=0.1$  and 
$$\eta_k=\frac{\|F_0\|}{(k+1)^2}.$$

Recall that the parameter $\theta_k$ is aimed at preventing $B_{k+1}$ from singularity. 
In all our numerical experiments, there was no single case, where this parameter differed from one. This means that all matrices $B_{k+1}$ generated by formula \eqref{QNupdates} were nonsingular. 

\begin{table}[thp]
\centering 	
\caption{List of test problems.}\label{t1}
\begin{tabular}[c]{lc}
	\hline
	\textsc{Problem} & \textsc{Dimension} \\
	\hline
	Brown almost-linear  & 10, \ 20, \ 30\\
	Broyden bounded & 10, \ 20, \ 30\\
	Broyden tridiagonal & 10, \ 20, \ 30\\
	Discrete boundary value & 10, \ 20, \ 30\\
	Discrete integral & 10, \ 20, \ 30 \\
	Trigonometric & 10, \ 20, \ 30\\
	Powell singulat & 4\\
	Helical valley & 3\\
	Powell badly scaled & 2\\
	Rosenbrock  & 2\\
	\hline
\end{tabular}
\end{table}

\begin{figure}[htp]
\vspace{-2.1in}
\includegraphics[width=1.00\textwidth]{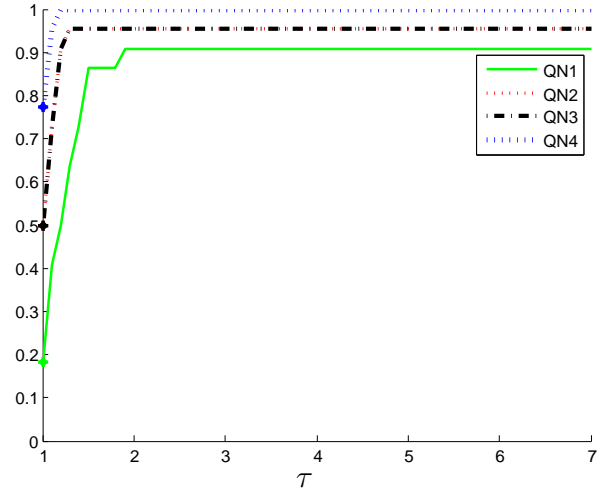} \\[-2.7in]
\begin{center} $\tau$ \end{center}
\caption{Performance profiles for the number of iterations. }
\label{fig_iter}
\end{figure}

\begin{figure}[htp]
\vspace{-2.9in}
\includegraphics[width=1.00\textwidth]{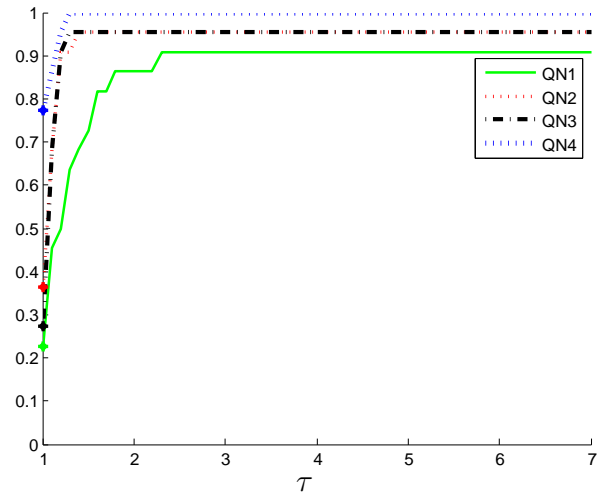} \\[-2.7in]
\begin{center} $\tau$ \end{center}
\caption{Performance profiles for the number of function evaluations. }
\label{fig_func_eval}
\end{figure}

For making experiments, we used 30 test problems from \cite{or-70}. They are listed in Table~\ref{t1}. The results of these experiments for the four algorithms are represented by the Dolan-Mor\'e performance profiles \cite{ds-96} based on the number of iterations, Fig.~\ref{fig_iter}, and the number of function evaluations, Fig.~\ref{fig_func_eval}. For $\tau = 1$, this performance measure indicates the portion of problems for which a given algorithm was the best. When $\tau > 1$, the profile, say, for the number of iterations, provides the portion of problems solved by a given algorithm in a number of iterations in each of these problems which does not exceed the $\tau$ times the number of iterations required by the algorithm that was the best in solving the same problem.

Recall that the computed values of $F(x)$ contains an information about the Jacobian matrix. Following the discussions in Section~\ref{sec:QN}, we sorted the algorithms from QN1 to QN4 in the way that they utilize this information more and more completely if to compare them in this order. The quality of the Jacobian approximation, which is related to the ability of reducing $\|F(x)\|$ along the corresponding search direction, improves following the suggested order of the algorithms. Figures~\ref{fig_iter} and \ref{fig_func_eval} illustrate how this quality affects the number of iterations and function evaluations. One can see that the best and worst performance was demonstrated by the interpolation method \cite{bf-94,b-97} and Broyden's method \cite{b-65}, respectively. The performance of the multipoint secant methods \cite{gs-78,b-83,b-86} was in between those associated with QN1 and QN4. Here it is necessary to draw attention to the robustness of the interpolation method.

As it was mentioned above, the multipoint and interpolation methods are mostly efficient in solving problems in which function evaluations are computationally more expensive than the linear algebra overheads associated with producing search directions. This is the reason why in our computer implementation of these methods we did not tend to reduce their CPU time. Therefore, we do not report here the time or running them. As expected, Broyden's method was the fastest in terms of time in 72\% of the test problems. However, it was less robust than the other methods.

\section{Conclusions}
\label{sec:Conclusions}
One of the main purposes of the present paper was to draw attention 
to the multipoint secant and interpolation methods as an alternative to Broyden's method. They were combined with the Li-Fukushima line search, and their global and superlinear convergence was proved. 

Our numerical experiments indicated that the multipoint secant and interpolation methods tend to be more robust and efficient than Broyden's method in terms of the number of iterations and function evaluations. This is explained by the fact that they are able to more completely utilize the information about the Jacobian matrix contained in the already calculated values of $F(x)$. It was observed that the more completely such information is utilized, the fewer iterations and number of function evaluations are, in general, required for solving problems. However, the linear algebra overheads related to the calculation of their search directions are obviously larger as compared with Broyden's method. Therefore, they can be recommended for solving problems with expensive function evaluations.

\section*{Acknowledgements}
Part of this work was done during Ahmad Kamandi's visit to Link\"oping University, Sweden. This visit was supported by Razi University.



\section*{References}
  \bibliographystyle{elsarticle-num} 
  \bibliography{references}





\end{document}